\documentclass[aos,MSNbibl,nameyear,seceqn,dvips]{arximspdf}
\usepackage{algorithmic,algorithm}

\doi{10.1214/14-AOS1228} 
\volume{42}
\issue{4}
\pubyear{2014}
\firstpage{1546}
\lastpage{1563}

\makeatletter
\newcommand{\new}{\mathrm{new}}

\newcommand{\Z}{{\mathbf Z}}
\newtheorem{theorem}{Theorem}[section]
\newproclaim{definition}[theorem]{Definition}
\newtheorem{lemma}[theorem]{Lemma}
\newtheorem{corollary}[theorem]{Corollary}
\newproclaim{example}[theorem]{Example}
\makeatother

\begin{document}
\begin{frontmatter}

\title{Optimum mixed level detecting arrays\thanksref{T1}}
\runtitle{Optimum mixed level detecting arrays}

\begin{aug}
\author{\fnms{Ce} \snm{Shi}\ead[label=e1]{shice060@163.com}\thanksref{t11}},
\author{\fnms{Yu} \snm{Tang}\corref{}\ead[label=e2]{ytang@suda.edu.cn}\thanksref{t22}}
\and
\author{\fnms{Jianxing} \snm{Yin}\ead[label=e3]{jxyin@suda.edu.cn}}
\runauthor{C. Shi, Y. Tang and J. Yin}
\affiliation{Soochow University}
\address{School of Mathematical Sciences\\
Soochow university\\
Suzhou, Jiangsu 215006\\
P.~R. China\\
\printead{e1}\\
\phantom{E-mail: }\printead*{e2}\\
\phantom{E-mail: }\printead*{e3}}
\end{aug}
\thankstext{T1}{Supported by NSFC Grants 11271280, 11271279,
11301342 and the NSF of Jiangsu Province Grant BK2012612.}
\thankstext{t11}{Supported by Shanghai Special Research Fund for
Training College's Young Teachers: ZZlx13001.}
\thankstext{t22}{Supported by Qing Lan Project.}

\received{\smonth{8} \syear{2013}}
\revised{\smonth{4} \syear{2014}}

%
\begin{abstract}
As a type of search design, a detecting array can be used to generate
test suites to identify and
detect faults caused by interactions of factors in
a component-based system. Recently, the construction and optimality of
detecting arrays have been investigated in depth in the case where all
the factors
are assumed to have the same number of levels.
However, for real world applications, it is more desirable to use
detecting arrays in which the various factors may have different
numbers of levels.
This paper gives a
general criterion to measure the optimality of a mixed level detecting array
in terms of its size. Based on this optimality criterion, the
combinatorial characteristics of mixed level detecting arrays of optimum
size are investigated. This enables us to construct
optimum mixed level detecting arrays
with a heuristic optimization algorithm
and combinatorial methods. As a result,
some existence results for optimum mixed level detecting arrays
achieving a lower bound
are provided for practical use.
\end{abstract}

%
\end{frontmatter}

\section{Introduction}\label{SecIntro}

Testing plays an important role in component-based systems. Due
to the complexity of systems, the number of possible
tests can be exponentially large.
Consider, for example,
the manufacturing of a personal computer system.
Suppose that a wide variety of choices are available for both hardware
and software components including the central unit processor,
primary memory, cache memory, interface,
operation system, web browser
and email system.
For manufacturers, if an exhaustive test suite is applied to
test whether there is any component interaction that may cause system
failure, the testing burden shall be extremely heavy even
when the system is of moderate complexity.
Consequently, they are interested
in generating test suites
that cover the most prevalent interactions.
This phenomenon is not unique and similar scenarios exist
in other disciplines such as
agriculture, manufacturing, networks and so on. For related
information,
readers are referred to \citet{CMMSSY2006},
\citet{MSSW2003}
and the references therein.

For certain component-based systems,
it is generally impossible
to examine every interaction without conducting
an exhaustive testing.
A popular and effective strategy is to set a limit
on the number of components involved in
interactions which may cause a failure.
Such an approach not only reduces the testing burden
dramatically, but
is also supported by many empirical results.
\citet{DKLPH}
claimed that a
large percentage of the existing faults in a software system
can be found only by examining all pairwise interactions,
while \citet{KWG} reported similar results.
For this reason, engineers and researchers prefer to
employ covering arrays
to conduct test suites.

A covering array (CA) is a special type of fractional factorial design.
The formal definition of a CA will be given in
Section~\ref{SecNot}.
Roughly speaking, a covering array of strength $t$ (index unity) is
a matrix, such that for every choice of $t$ columns of the matrix,
all possible $t$-tuples appear within its rows.
When applying a CA to testing problems, the columns of the CA
are utilized to represent factors affecting a response in which
the entries within the columns indicate settings or values for that
factor. The rows represent tests to be performed, in which a value
for each factor is dictated.
Covering arrays have been extensively studied;
see \citet{Col2004} for a survey for CAs,
and \citet{CMMSSY2006},
\citet{MSSW2003} and \citet{CSWY2011} for their properties,
applications and efficient
constructions.

Testing with a CA can determine
whether any
$t$-component interactions cause a failure,
which is an important step in screening a system for interaction faults.
However, as \citet{CM2008} pointed out,
test suites based on covering arrays provide
little information about identifying the exact interactions causing
faults. Due to practical
concerns, tests that reveal the location of interaction faults are of
interest. \citet{CM2008} formalized the problem of nonadaptive
location of interaction faults under the hypothesis that the system
contains (at most) $d$ faults, each involving (at
most) $t$ factors.
They proposed the notion of a detecting array, which is a CA with
additional requirements.

Under the framework of \citet{CM2008},
detecting arrays have been investigated
in depth by a number of authors in the case where
all the factors had the
same number of levels. A combinatorial characterization
and a general criterion of
measuring optimality for such detecting arrays were
established in \citet{STY2011}
and \citet{TY2011}. For practical
applications, things are not always that simple, and it is
desirable to use detecting arrays in which the various factors
may have different
numbers of levels.

In \citet{BPP1992}, the authors considered a test
on the copy function of the PMx/StarMAIL system,
aiming to provide a pairwise covering solution.
Four factors, that is, function scope, server type, client type and
target content
were considered, each with three levels.
The factor ``function scope'' referred to three copy
options: copy only a single (current) message;
copy messages previously marked, or copy all messages contained
within a folder. However, in practice it is easy to use
some shortcut keys to mark all messages in a folder, thus
we can only assign two values, that is, ``current'' and ``marked,''\vadjust{\goodbreak} to
the factor ``function scope,'' as done by many other electronic mail systems.
In view of this,
the original fixed level test problem
becomes a mixed level one.
Another mixed level test example can be found in
\citet{CDPP1996}, where the authors considered a voice response
unit test,
affected by four factors with three, three, two and two levels.

In the above two examples, the authors both
suggested employing an orthogonal array or
a covering array to conduct the corresponding test suites.
As a result, interaction faults in these systems
can be found,
but their locations cannot always be identified.
Such drawbacks can be
overcome by using a detecting array.
Table~\ref{TabFirstEx} shows an example of a detecting array
for testing the copy function of the PMx/StarMAIL system
(here ``function scope'' is assigned to be a two-level factor).
The array contains 18 tests.
If we assume that an error in this system can only be due to
some combination of two components, then any single fault can
be identified according to the outcome of these tests.
For example, the last column of Table~\ref{TabFirstEx}
lists a possible outcome, which indicates the system error must be
caused by copying the current message by an MSNET client.
Moreover, if there is more than one combination of two components
causing the faults, this can also be detected.

\begin{table}
\tabcolsep=0pt
\caption{A detecting array for the copy function of the PMx/StarMAIL system}\label{TabFirstEx}
\begin{tabular*}{\tablewidth}{@{\extracolsep{\fill}}@{}lccccc@{}}
\hline
\textbf{Test} & \textbf{Function scope} & \textbf{Server type} & \textbf{Client type} & \textbf{Target content} & \textbf{Outcome}\\
\hline
\phantom{0}1 & {Current} & {Share mode} & {MSNET} & {Empty} & {F}\\
\phantom{0}2 & {Current} & {User mode} & {MSNET} & {Partial} & {F}\\
\phantom{0}3 & {Marked} & {Share mode} & {Enhanced} & {Partial} & {P}\\
\phantom{0}4 & {Marked} & {User mode} & {Basic} & {Partial} & {P}\\
\phantom{0}5 & {Marked} & {MSNET} & {MSNET} & {Partial} &{P}\\
\phantom{0}6 & {Marked} & {MSNET} & {Basic} & {Full} &{P}\\
\phantom{0}7 & {Current} & {Share mode} & {Enhanced} & {Full} & {P}\\
\phantom{0}8 & {Marked} & {User mode} & {MSNET} & {Empty}& {P}\\
\phantom{0}9 & {Current} & {User mode} & {Enhanced} & {Empty} & {P}\\
10 & {Marked} & {MSNET} & {Enhanced} & {Empty}& {P}\\
11 & {Current} & {Share mode} & {Basic} & {Partial} & {P}\\
12 & {Current} & {MSNET} & {Basic} & {Empty} &{P}\\
13 & {Marked} & {Share mode} & {MSNET} & {Full}& {P}\\
14 & {Marked} & {Share mode} & {Basic} & {Empty} & {P}\\
15 & {Current} & {MSNET} & {Enhanced} & {Partial} & {P}\\
16 & {Current} & {MSNET} & {MSNET} & {Full} &{F}\\
17 & {Marked} & {User mode} & {Enhanced} & {Full} & {P}\\
18 & {Current} & {User mode} & {Basic} & {Full}& {P}\\
\hline
\end{tabular*}
\end{table}

If properly coded, a detecting array is actually
a search design without noise [\citet{Srivastava1975}]. Let $S$ be the
set of all possible combinations of
certain components causing faults (for the detecting array in Table~\ref
{TabFirstEx}, $S$ is the
set of all { $3\times9+ 3 \times6 = 45$} level combinations with two
components),
$Y$ represent the outcome, that is, a 0--1 vector (zero means failure,
while one means pass),
and $M$ be the incidence matrix of the array (denoted by $A$) related
to $S$,
that is, if a level combination in $S$ occurs in a row of $A$, the
corresponding entry takes the value one, otherwise it takes the value zero
(in the above example, $M$ shall be an {$18$ by $45$} zero--one matrix).
If $A$ is a detecting array, then the solution for $\beta$ in equation
$M \beta= Y$
is unique over a predefined space for any possible outcome $Y$
(for the detecting array in Table~\ref{TabFirstEx}, the solution
space for $\beta$ only contains
45 unit vectors with one element 1 and others 0). Notice here
the special coding of $M$ and the restricted space of $\beta$
make it different from the usual search designs discussed
in the existing literature;
see, for example,
\citet{STS1996}, \citet{GB2001} and \citet{ETMJ2011}.
Under this consideration, the present authors took a new approach
to investigate the property and requirement for detecting arrays.
In \citet{TY2011} and \citet{STY2011},
they obtained many interesting results
for {fixed level detecting arrays}.
Although the basic concepts and notation related to
fixed level detecting arrays are analogous
to those of mixed level detecting arrays,
the combinatorial characteristics as well as construction methods
cannot be directly duplicated
due to the complicated structure of mixed level arrays.
In fact, optimum mixed level detecting arrays are not
always
super-simple (see Definition \ref{SuperS}),
but have the ``extendible'' property,
which will be defined in Section~\ref{SecOpt}.
Based on this special property, we will further provide
a heuristic algorithm and some combinatorial methods
to construct optimum mixed level arrays with specific parameters.

The organization of this paper
is as follows. In Section~\ref{SecNot},
we give a detailed description
on the definition of
a mixed level detecting array. In Section~\ref{SecOpt}, we establish a
general criterion for measuring the optimality of a mixed
level detecting array in terms of its size, and then describe the
combinatorial characteristics of mixed level detecting arrays of optimum
size. Based on properties in Section~\ref{SecOpt}, a heuristic optimization algorithm
to generate mixed level detecting arrays is developed
in Section~\ref{sec4} and some optimum detecting arrays with the
largest number of factors
are then found. Composition methods for constructing
optimum detecting arrays based on combinatorial methods
are included in Section~\ref{sec5}. This provides a number of infinite classes of
optimum mixed level detecting arrays.
Section~\ref{s6} concludes, while all proofs are deferred to the \hyperref[app]{Appendix}.

\section{Definitions and terminology}\label{SecNot}

Let $v_1, v_2, \ldots, v_{k}$ be $k$ natural numbers
(not necessarily distinct). For each $i$ with $1\le i\le k$, let $V_i$
be a set of cardinality $v_i$. For given natural numbers $N$ and $t\le k$,
a \emph{mixed covering array} (MCA) of type $(v_1,v_2,\ldots, v_k)$,
size $N$,
strength $t$ and index $\lambda$, denoted by $\mathrm{MCA}_{\lambda}(N; t, k,
(v_1,v_2,\ldots, v_k))$,
is an $N\times k$ array $A$, which satisfies the following two properties:
\begin{longlist}[(2)]
\item[(1)] For each $i$ with $1\le i\le k$, the entries in the $i$th column
of $A$ are taken from~$V_i$.
\item[(2)] The rows of each $N \times t$ subarray of $A$ cover all $t$-tuples
of values from the $t$ columns at least $\lambda$ times.
\end{longlist}
The coordinates of $V_1 \times V_2 \times\cdots\times V_k$
are referred to as \emph{factors},
so the elements of $V_i$ represent the \emph{levels} of factor $i$
for $1\le i\le k$ and the parameter $k$ is the number of factors, also called
the \emph{degree}. The mixed covering array
number MCAN$_\lambda(t,k,(v_1,v_2,\ldots, v_k))$ is the minimum $N$ required
to produce an $\mathrm{MCA}_\lambda(N; t,k,(v_1,v_2,\ldots, v_k))$. A mixed
covering array
is termed \emph{optimum} if its size $N=\mathrm{MCAN}_\lambda
(t,k,(v_1,v_2,\ldots, v_k))$.

The term ``mixed'' in the definition is used to indicate that the $k$
level numbers $v_i$ $(1\le i\le k)$
may take different values.
If the MCA has a fixed number of levels, it is often known as a
\emph{covering array}. In this case, the notation CA$_{\lambda}(N; t,
k, v)$
and CAN$_\lambda(t,k, v)$ is employed. In the literature, the
subscript is often
dropped from the above notation whenever $\lambda=1$.
A~CA$_{\lambda}(N; t, k, v)$ of size $N=\lambda v^t$ is known as an
\emph{orthogonal array} (OA), or an $\mathrm{OA}_{\lambda}(t,k,v)$, if the rows
of every $N\times t$
subarray cover all $t$-tuples of symbols exactly $\lambda$ times. Clearly,
an $\mathrm{OA}_{\lambda}(t,k,v)$
is an {optimum} CA$_{\lambda}(N; t, k, v)$ of size $N = \lambda v^t$.

Notice that the property of an $\mathrm{MCA}_\lambda(N; t,k,(v_1,v_2,\ldots, v_k))$
is preserved if
its columns are permuted.
So, the type $(v_1,v_2,\ldots, v_k)$ of
an MCA is essentially an unordered multiset and
one can assume
that the sizes of the $k$ level sets are in a nondecreasing order, that is,
$v_1\le v_2\le\cdots\leq v_k$.
It follows that an $\mathrm{MCA}_\lambda(N; t,k,(v_1,v_2,\ldots, v_k))$ can
exist only if
$N\geq\lambda\prod_{i=k-t+1}^k v_i$. It is {optimum} when
$N=\lambda\prod_{i=k-t+1}^k v_i$.
For convenience, a shorthand
notation is adopted to describe the type $(v_1,v_2,\ldots, v_k)$
by combining $v_i$'s that are the same. For example,
if three $v_i$'s are of size two, one writes this as $2^3$.

We use $I_k$ to denote the set of the first $k$ natural numbers.
In order to describe the definition of a detecting array,
let us formulate the testing problem formally.
Consider a system with $k$ factors (parameters or components).
One may identify the set of
values (levels) of the $j$th
factor $(1\le j\le k)$ with $\Z_{v_j}=\{0,1, \ldots, v_j-1\}$,
the residue class ring of integers modulo $v_j$.
Suppose that $v_1 \le v_2 \le\cdots\le v_k$.
A~test is an assignment of values to factors, that is,
a $k$-tuple from $\Z_{v_1}\times\Z_{v_2}\times \cdots\times\Z_{v_k}$. The
execution of a test can have two outcomes: \emph{pass} or \emph{fail}. A $t$-way interaction or an interaction of strength $t$ is a set
of the form $\{(j_i,\sigma_{j_i})|1\le i
\le t\}$, where $\{j_1, j_2, \ldots, j_t\} \subseteq I_k$ with
$j_1<j_2 < \cdots< j_t$ and $\sigma_{j_i} \in\Z_{v_{j_i}}$.
A test (or a \mbox{$k$-}tuple) $R= (x_1, \ldots, x_k) \in
\Z_{v_1}\times\Z_{v_2}\times\cdots\times\Z_{v_k}$ covers the interaction
$\{(j_i,\sigma_{j_i})|1\le i
\le t\}$ if $x_{j_i} = \sigma_{j_i}$ for $ 1 \leq i \leq t$.
Thus, a test with $k$ factors covers exactly ${k\choose t}$
interactions of strength $t$. A test suite is a collection of
tests; the outcomes are the corresponding set of pass/fail results.
We assume that failures are caused by faulty interactions. A~test
fails when it contains at least one of the faulty
interactions, and does not fail otherwise. It then turns out that
a test suite is essentially an $N\times k$ array $A$ of type
$(v_1,v_2,\ldots, v_k)$ whose rows consist
of $N$ elements (possibly with multiplicities) of the Cartesian product
$\Z_{v_1}\times\Z_{v_2}\times\cdots\times\Z_{v_k}$.
In order to observe an
interaction fault, it is necessary that the interaction is covered
by at least one test in the test suite. This is to say that the array $A$
is required to be
an $\mathrm{MCA}_{\lambda}(N; t, k, (v_1,v_2,\ldots, v_k))$.

Now suppose that $A=(a_{ij})$ $(i \in
I_N, j\in I_k)$ is an $\mathrm{MCA}_{\lambda}(N; t, k, (v_1, v_2, \ldots,\break  v_k))$
whose entries on the $j$th column are taken from $\Z_{v_j}$ for $1\le
j\le k$.
As we have noted before, $t$-way interactions that cause faults
cannot be located by testing with $A$ as suite.
Thus, in order to determine and identify
the faults, the covering array $A$ must possess
more structure than stipulated in its
definition. To this end, we define $\rho(A, T)$ to be
the set of rows of $A$ in each of which, an arbitrary $t$-way interaction
$T$ is covered.
To be more precise, for any $t$-way interaction $T = \{(j_i,\sigma
_{j_i})|1\le i
\le t\}$, define
\[
\rho(A, T) = \{r| a_{rj_i}=\sigma_{j_i}, 1\le i\le t\},
\]
where $a_{rj_i}$ is the $(r,j_i)$ entry of $A$. For a set of
interactions ${\mathcal T}$, define
\[
\rho(A, {\mathcal T}) = {\bigcup_{T \in{\mathcal T}} \rho(A, T)}.
\]
Following \citet{CM2008},
an $\mathrm{MCA}_{\lambda}(N; t, k, (v_1,v_2,\ldots, v_k))$ $A$ is termed
a detecting array (DTA), or
a $(d,t)$-$\mathrm{DTA}(N; k, (v_1,v_2, \ldots, v_k))$, if
%
\begin{equation}
\label{DTAdef} \rho(A, T) \subseteq\rho(A, {\mathcal T})
\quad\Longleftrightarrow\quad T \in{
\mathcal T},
\end{equation}
whenever $T$ is a $t$-way interaction and ${\mathcal T}$ is a set of
$t$-way interactions of cardinality~$d$.

It was proved in \citet{CM2008}
that testing with a $(d,t)$-$\mathrm{DTA}(N; k, (v_1,v_2, \ldots, v_k))$
is able to identify any set of $d$
interaction-faults of strength $t$ from the
outcomes. Further, if there are more than $d$ $t$-way interactions
causing the faults, {this} can also be detected.

As with MCAs, the minimum number $N$ of runs for which a
$(d,t)$-$\mathrm{DTA}(N; k, (v_1,v_2,\ldots, v_k))$ exists is called \emph{detecting array number}, denoted by $(d,t)$-$\mathrm{DTAN}(k, (v_1, v_2,
\ldots, v_k))$. A $(d,t)$-$\mathrm{DTA}(N; k, (v_1, v_2, \ldots, v_k))$ with
$N=(d,t)$-$\mathrm{DTAN}( k, (v_1, v_2, \ldots, v_k))$ is said to be \emph{optimum}. In the case $v_1=v_2=\cdots
= v_k=v$, the notation $(d,t)$-$\mathrm{DTA}(N; k, v)$ and $(d,t)$-$\mathrm{DTAN}(k, v)$
are adopted.

Remark that the detecting array of Table~\ref{TabFirstEx} in Section~\ref{SecIntro}
is a
$(1,2)$-$\mathrm{DTA}(18; 4, (2^13^3))$ if its levels are
mapped to the corresponding residue class rings of integers.
It can be readily checked that
for any two distinct 2-way interactions $T$ and $T'$,
we have $\rho(A, T) \ne\rho(A, T')$. This means that condition
(\ref{DTAdef})
is satisfied, since $\rho(A, T) \ne\rho(A, T')\Longrightarrow T \ne T'$.
The detecting array is of strength $t=2$.
All its 2-way interactions among columns are covered
at least twice.
So it is actually an $\mathrm{MCA}_{2}(18; 2, 4, (2^13^3))$.
Furthermore, the rows of any $t+1=3$ columns cover every triple
of values from the $3$ columns at most once. This is a
feature of the underlying MCA of a DTA.
To characterize this feature, we introduce the following concept,
which will be frequently used in the subsequent sections.

\begin{definition}\label{SuperS}
An $\mathrm{MCA}_\lambda( N;t,k, (v_1, v_2, \ldots, v_k))$ is said to be \emph{super-simple},
if each of its $N\times(t+1)$ subarray covers every $(t+1)$-tuple of
values from the $t+1$ columns at most once.
\end{definition}

\section{{Optimality criterion and combinatorial feature}}\label{SecOpt}

This section is devoted to establishing a general criterion
for measuring the optimality of DTAs, and giving the
combinatorial characterization of mixed level DTAs with {optimum} size.
Suppose that $A$ is an arbitrary $(d,t)$-$\mathrm{DTA}(N; k, (v_1, v_2,\break  \ldots,v_k))$ with $v_1\le v_2\le\cdots\leq v_k$.
We shall exclude the following three trivial cases.
The first one is $t=k$.
In this case, an
optimum DTA
is just the full factorial design.
The second case is that of $v_1=1$ for which the DTA is determined uniquely
by its subarray, the one with the first
factor removed. The last case is $d \ge v_1$.
In this case, \citet{CM2008}
showed that a $(d,t)$-$\mathrm{DTA}(N; k, (v_1, v_2, \ldots,
v_k))$ cannot exist.
So the restrictions
\[
2\le v_1\le v_2\le\cdots\leq v_k,\qquad t < k\mbox{ and } d < v_1
\]
are always assumed, whenever
we speak of a $(d,t)$-$\mathrm{DTA}(N; k, (v_1, v_2, \ldots,v_k))$. In addition, the $j$th level-set
of cardinality $v_j$ is understood to be $\Z_{v_j}$ for $1\le j\le k$,
unless otherwise stated.
Under these conventions,
we have the following lemma which is an analogue of Theorem 2.3 in \citet{STY2011}.
Hence, we omit its proof for brevity.

\begin{lemma} \label{lbound-d>1}
Suppose that $A$ is a $(d,t)$-$\mathrm{DTA}(N;k,(v_1, v_2, \ldots, v_k))$.
Then $|\rho(A,T)|\geq d+1 $ for any $t$-way interaction $T$.
\end{lemma}

By invoking Lemma \ref{lbound-d>1}, we are able to establish a lower bound
on the function $(d,t)$-$\mathrm{DTAN}( k, (v_1, v_2,\ldots, v_k))$, which
serves as
a general criterion
for measuring the optimality of detecting arrays.

\begin{theorem}\label{Thm-LowB-M}
Let $v_j$ $(1\le j \le k)$ be $k$ integers with $2\leq v_1\leq v_2\leq
\cdots\leq v_k$. Then
%
\begin{equation}
\label{Eqn-LowB-M} (d,t)\mbox{-}\mathrm{DTAN}\bigl(k, (v_1, v_2,
\ldots, v_k)\bigr) \ge(d+1)\prod_{i=k-t+1}^k
v_i.
\end{equation}
\end{theorem}

Taking $v_1=v_2=\cdots=v_k=v$ in Theorem \ref{Thm-LowB-M}, we
retrieve the lower
bound on the function $(d,t)$-$\mathrm{DTAN}(k,v)$ given in
\citet{STY2011}.

\begin{corollary}\label{Thm-LowB-F}
Let $t, k$ and $v$ be positive integers with $t < k$. Then
%
\begin{equation}
\label{Eqn-LowB-F} (d,t) \mbox{-}\mathrm{DTAN} (k,v) \geq(d+1)v^t.
\end{equation}
\end{corollary}

It was proved in \citet{STY2011} that an optimum $(d,t)$-\break $\mathrm{DTAN}(N; k,v)$
meeting the lower bound in (\ref{Eqn-LowB-F}) is
equivalent to a super-simple $\mathrm{OA}_{d+1}(t, k,v)$. However, the structure
of an {optimum}, mixed level DTA is much more complicated
than that of an {optimum}, fixed level DTA
even though the lower bounds in (\ref{Eqn-LowB-M})
and (\ref{Eqn-LowB-F}) look similar.
The aforementioned equivalence does not hold
in the mixed level case, which can be seen in the following example.

\begin{example}\label{ex-ODTA2}
An optimum $(2,2)$-$\mathrm{DTA}(48; 5, (3^34^2))$ is formed
by taking the transpose of the superimposition of the following two
$24\times5$ arrays:
{\fontsize{9.5pt}{11pt}\selectfont{\begin{eqnarray*}
\pmatrix{1 & 2 & 0 & 0 & 2 & 0 & 2 & 2 & 1 & 0 & 0 &1 & 1 & 2& 1 & 2 & 0 & 2 & 0 & 0 & 1 & 0 & 1 & 2\cr
2 & 1 & 2 & 0 & 0 & 2 & 2 & 1 & 0 & 1 & 1 & 2 & 0 & 1 & 1 & 2 & 2 & 2 &0 & 1 &0 & 2 & 1 & 2\cr
1 & 0 & 0 & 2 & 1 & 2 & 2 & 1 & 2 & 1 & 1 & 1 & 1 & 2 & 0 & 2 & 0 & 0 &1 & 0 &0 & 0 & 0 & 1\cr
2 & 3 & 1 & 2 & 1 & 3 & 3 & 0 & 1 & 1 & 2 & 3 & 0 & 2 & 0 & 1 & 3 & 0 &2 & 1 &3 & 2 & 0 & 0\cr
1 & 2 & 1 & 1 & 1 & 0 & 1 & 2 & 0 & 2 & 3 & 2 & 3 & 0 & 1 & 2 & 3 & 3 &0 & 0 &1 & 2 & 0 & 0},
\\[3pt]
\pmatrix{0 & 0 &2 & 1 & 2 & 2 & 1 & 2 & 1 & 1 & 2 & 0 & 1 & 2 & 2 & 0 & 1 & 1 & 1 & 0 &0 & 0 & 2& 1\cr
1 & 0 & 0 & 2 & 2 & 1 & 1 & 0 & 0 & 1 & 1 & 0 & 2 & 0 & 2 & 1 & 0 & 1 &2 & 0 &0 & 2 & 0 & 1\cr
1 & 2 & 2 & 0 & 2 & 0 & 2 & 1 & 0 & 1 & 0 & 2 & 2 & 1 & 1 & 2 & 0 & 2 &1 & 2 &0 & 1 & 0 & 2\cr
3 & 0 & 0 & 2 & 2 & 1 & 2 & 2 & 1 & 3 & 2 & 3 & 0 & 3 & 1 & 0 & 2 & 3 &1 & 1 &0 & 0 & 3 & 1\cr
1 & 0 & 1 & 0 & 3 & 3 & 2 & 2 & 2 & 0 & 1 & 2 & 2 & 3 & 0 & 3 & 3 & 3 &3 & 3 &2 & 1 & 0 & 1}.
\end{eqnarray*}}}%
It is readily checked that this detecting array is not super-simple.
\end{example}

Since an optimum, mixed level DTA
may not be super-simple,
we define more set-theoretic
configurations to explore
the combinatorial features
of an optimum, mixed level DTA.

\begin{definition}\label{DefExtension}
Let $T_{t}=\{(j_i,\sigma_{j_i})|1\le i \le t\}$ be a
$t$-way interaction. Then each $(t+1)$-way interaction
$T_{t+1}=((j_1,\sigma_{j_1}),\ldots,(j_t,\sigma_{j_t}),
(j_{t+1},\sigma_{j_{t+1}}))$
is said to be an extension of $T_{t}$, where $ j_{t+1}\in I_k
\setminus
\{j_1,j_2,\ldots, j_t\}$.
\end{definition}

The following lemma follows Definition
\ref{DefExtension} immediately.

\begin{lemma} \label{T,T+1}
Suppose that $A$ is an $N\times k$ array of type $(v_1, v_2, \ldots, v_k)$.
Let $T_t$ be a $t$-way interaction and
$T_{t+1}$ an extension of $T_t$. Then $\rho(A,T_{t+1})
\subseteq\rho(A,T_t)$.
\end{lemma}

Lemma \ref{T,T+1} motivates us
to introduce the following concept.

\begin{definition}\label{Defextend}
Let $A$ be an $N\times k$ array of type $(v_1, v_2, \ldots, v_k)$.
If for any $t$-way interaction $T_{t}$ of $A$ and
any $d$ {extension} $T_{t+1,i}$ $(1 \le i \le d)$
of $T_t$, we have
$\rho(A,T_{t}) \setminus\bigcup_{i=1} ^d
\rho(A,T_{t+1,i}) \neq\varnothing$, then the array
is said to be \emph{$d$-{extendible}}.
\end{definition}

It is interesting to note that in the particular case where the array $A$
is an $\mathrm{OA}_{d+1}(t,k,v)$, the terms ``$d$-{extendible}'' and ``super-simple''
are essentially the same. Specifically, we have the following lemma.

\begin{lemma}\label{TheOAEqu}
Let $A$ be an $\mathrm{OA}_{d+1}(t,k,v)$. Then $A$ is $d$-extendible
if and only if $A$ is super-simple.
\end{lemma}

Now we are in a position to characterize the combinatorial features
of an optimum $(d,t)$-$\mathrm{DTA}(N; k, (v_1, v_2, \ldots, v_k))$.

\begin{theorem}\label{TheEquivd}
Suppose that $A$ is an $N\times k$ array of type $(v_1, v_2, \ldots, v_k)$.
Then $A$ is a $(d,t)$-$\mathrm{DTA}(N; k, (v_1, v_2, \ldots, v_k))$
if and only if $A$ is a $d$-extendible $\mathrm{MCA}_{d+1}(N;t, k, (v_1, v_2,
\ldots,
v_k))$.
\end{theorem}

As an immediate consequence of Theorem \ref{TheEquivd}, we have the
following result.

\begin{theorem}\label{TheOptimalEquivd}
An optimum $(d,t)$-$\mathrm{DTA}(N; k, (v_1, v_2, \ldots, v_k))$
meeting the lower bound in (\ref{Eqn-LowB-M}) is
equivalent to a $d$-{extendible} $\mathrm{MCA}_{d+1}(N;t, k, (v_1,\break  v_2,
\ldots,
v_k))$ of optimum size $N=(d+1)\prod_{i=k-t+1}^k v_i$.
\end{theorem}

Combining Lemma \ref{TheOAEqu} with Theorem \ref{TheOptimalEquivd}
gives us the following useful corollary, which was first stated in \citet{STY2011}.

\begin{corollary}\label{fixedlevel}
An optimum $(d, t)$-$\mathrm{DTA}((d + 1)v^t; k, v)$
meeting the lower bound in (\ref{Eqn-LowB-F}) is
equivalent to
a super-simple $\mathrm{OA}_{d+1}(t, k, v)$.
\end{corollary}

The following two theorems take care of the cases $k=t+1$ and the
type $(v_1, v_2, \ldots, v_{t+1})= (a^1b^{k-1})$, respectively.\vadjust{\goodbreak}

\begin{theorem}\label{Efeature1}
A $(d,t)$-$\mathrm{DTA}(N; t+1, (v_1, v_2, \ldots, v_{t+1}))$ of optimum size $N$
meeting the lower bound in (\ref{Eqn-LowB-M}) is
equivalent to a super-simple
$\mathrm{MCA}_{d+1}(N;\break  t, t+1, (v_1, v_2, \ldots, v_{t+1}))$ of {optimum}
size $N=(d+1)\prod_{i=2}^{t+1} v_i$.
\end{theorem}

\begin{theorem}\label{Efeature2}
Let\vspace*{2pt} $a$ and $b$ be positive integers satisfying $a\le b$. Then,
an optimum $(\lambda-1,t)$-$\mathrm{DTA}(N; k,(a^1b^{k-1}))$ meeting the
lower bound in (\ref{Eqn-LowB-M}) is
equivalent to a super-simple $\mathrm{MCA}_\lambda(N;t, k, (a^1b^{k-1}))$
of optimum size\break  $N=\lambda b^t$.
\end{theorem}

Finally, we give a sufficient condition for the existence of an
optimum detecting array, meeting the
lower bound in (\ref{Eqn-LowB-M}).

\begin{theorem}\label{Feature}
If a super-simple $\mathrm{MCA}_{d+1}(N;t, k, (v_1, v_2, \ldots,
v_k))$ of optimum size $N=(d+1) \prod_{i=k-t+1}^k v_i$ exists, then an
optimum $(d,t)$-$\mathrm{DTA}(N;\break  k, (v_1, v_2,
\ldots, v_k))$ meeting the lower bound in (\ref{Eqn-LowB-M}) also exists.
\end{theorem}

\section{A heuristic algorithm}\label{sec4}
In this section, we present an algorithm to search for
optimum detecting arrays meeting the lower bound in (\ref{Eqn-LowB-M})
with a small number of runs.
The basic idea is to use a heuristic optimization algorithm,
which is based on the combinatorial features described in the previous
section. We first
provide two necessary conditions for the existence of an
optimum detecting array, meeting the lower bound in (\ref{Eqn-LowB-M}).

\begin{lemma}\label{Lem-Fixed}
If there exists an optimum $(1,2)$-$\mathrm{DTA}(2q^2; k, q)$ meeting the lower
bound in (\ref{Eqn-LowB-M}),
then $k \le2q$.
\end{lemma}

\begin{lemma}\label{Lem-Existence}
Let $u, k$ and $w$ be nonnegative integers and $w \ge3$.
If there exists an optimum $(1,2)$-$\mathrm{DTA}(N; 1+k+u, 2^u 3^k w^1)$
meeting the lower bound in (\ref{Eqn-LowB-M}),
then either $u$ or $k$ is less than $2$. Moreover,
if $k=0$, $u \le3$; if $k=1$, $u \le4$;
if $u=0$, $k \le5$; if $u=1$, $k \le3$.
\end{lemma}

For parameters $N$, $k$, $v_1, v_2, \ldots, v_k$
satisfying Lemmas \ref{Lem-Fixed} or \ref{Lem-Existence} constraints,
we randomly generate an $N \times k$ matrix, denoted by $A_0$, with
each of its $j$th column taking values
from $Z_{v_j}$. For any 2-way interaction $T$ and any of its possible
extension $T_a$, record the
values of $|\rho(T,A_0)|$ and $|\rho(T_a,A_0)|$.
Given that $A_0$ is a $(1,2)$-detecting array if and only if
these two values are not equal, we then define the quantity
$\sum_{T,T_a} I_A (|\rho(T,A_0)| = |\rho(T_a,A_0)| )$ as
an objective function to minimize,
where $I_A(\cdot)$ represents the indicator function. We use the
simulated annealing method to find
a final design with the objective function taking the value zero.
The pseudo code in Algorithm~\ref{Alg1} illustrates the details.

\begin{algorithm}[t]
\begin{algorithmic}
\caption{Pseudo code for searching $(1,2)$-detecting arrays}\label{Alg1}
{\fontsize{10.5pt}{11pt}\selectfont{
\STATE Randomly generate $A_0$ and calculate {$\Delta(A_0) =\sum
_{T,T_a} I_A (|\rho(T,A_0)| = |\rho(T_a,A_0)| )$}
\WHILE{$\Delta(A_0) > 0$}
\STATE Randomly exchange two entries in the same column of $A_0$ to get
$A_{\new}$
\STATE Compute $\nabla= \Delta(A_{\new}) - \Delta(A_0)$
\IF{$\nabla< 0$}
\STATE Set $A_0 = A_{\new}$
\ELSE
\STATE Randomly generate $u$ from uniform distribution $U[0,1]$
\IF{$u < p = 0.01$}
\STATE Set $A_0 = A_{\new}$
\ENDIF
\ENDIF
\ENDWHILE}}%
\end{algorithmic}
\end{algorithm}

\begin{table}[b]
\tabcolsep=0pt
\caption{Optimum $(1,2)$-$\mathrm{DTA}(N;k,(v_1,v_2,\ldots,v_k))$ for $N\leq 30$}\label{TabResList}
\begin{tabular*}{\tablewidth}{@{\extracolsep{\fill}}@{}lcccccccc@{}}
\hline
$\bolds{N}$ & \textbf{8} & \textbf{12} & \textbf{16} & \textbf{18} & \textbf{18} & \textbf{18} & \textbf{20} & \textbf{24} \\
\hline
\mbox{Levels} & $2^a$ & $2^a3^1$ & $2^a4^1$ & $2^13^a$ & $2^a3^2$ &$3^a$ & $2^a5^1$ & $2^a6^1$ \\
\mbox{Existence} & $a\leq4$ & $a\leq3$ & $a\leq3$ & $a\leq3$ &$a\leq4$ & $a\leq6$ & $a\leq3$ & $a\leq3$
\\[6pt]
$\bolds{N}$ & \textbf{24} & \textbf{24} & \textbf{24} & \textbf{28} & \textbf{30} & \textbf{30} & \textbf{30} & \\\hline
\mbox{Levels} & $2^a3^24^1$ & $2^a3^14^1$ & $3^a4^1$ & $2^a7^1$ &$2^a3^35^1$ & $2^a3^15^1$ & $3^a5^1$ & \\
\mbox{Existence} & $a\leq1$ & $a\leq4$ & $a\leq5$ & $a\leq3$ &$a\leq1$ & $a\leq4$ & $a\leq5$ & \\\hline
\end{tabular*}
\end{table}

We list our search results for $N \le30$
in Table~\ref{TabResList}.
All the obtained DTAs are of optimum size meeting
the bound in (\ref{Eqn-LowB-M}).
These DTAs are not only useful in
real world applications but also
as ingredients in the construction of new optimum DTAs
via combinatorial approaches developed in the next section. The found
$(1,2)$-$\mathrm{DTA}(18; 4, (2^13^3))$
indicates the optimality
of the test scheme presented in Section~\ref{SecIntro}
regarding the copy function of PMx/StarMAIL
system. Meanwhile, the $(1,2)$-$\mathrm{DTA}(18;4,(2^2 3^2))$
in Table~\ref{TabResList} gives an optimum solution
to the voice response unit test problem (without additional restrictions)
studied by \citet{CDPP1996}.

We note that an exhaustive computer search shows
that neither a $(1,2)$-$\mathrm{DTA}(18;5,(2^1 3^4))$
nor a $(1,2)$-$\mathrm{DTA}(24;5,(2^1 3^3 4^1))$
can exist, although their \mbox{parameters} are compatible with Lemma \ref
{Lem-Existence} constraints.

\section{Combinatorial approaches}\label{sec5}

In real world applications,
some testing problems may have factors
with a large number of levels. For example,
suppose a program consists of several functions or subroutines.
Each function needs an input parameter and each parameter may take
hundreds of values. This
causes the number of factor level combinations to be very large.
The objective of this section is to present some combinatorial
approaches to find
optimum, mixed level detecting arrays with a large number of levels
per factor for practical applications.

\subsection{New optimum DTAs from the Kronecker product}

For given two matrices $A=(a_{ij})$ $(i \in I_{m}, j\in I_k)$ and
$B=(b_{rs})$ $(r \in I_{n},
s\in I_k)$, the
Kronecker product of $A$ and $B$ is defined to be the $(mn)\times k$ matrix
\begin{eqnarray*}
A\otimes B= \pmatrix{
(a_{11}& a_{12}& \cdots& a_{1k})\otimes B\cr
(a_{21}& a_{22}& \cdots& a_{2k})\otimes B\cr
\vdots& \vdots& \vdots& \hspace*{6pt}\vdots\hspace*{13pt}\vdots\hspace*{8pt}\vdots\cr
(a_{m1}& a_{m2}& \cdots& a_{mk})\otimes B},
\end{eqnarray*}
where for any $i$ with $1\le i\le m$, $(a_{i1}, a_{i2}, \ldots,a_{ik})
\otimes B$ is
defined to be the matrix
\begin{eqnarray*}
\pmatrix{(a_{i1},b_{11})& (a_{i2}, b_{12})& \cdots &(a_{ik}, b_{1k})\cr
(a_{i1}, b_{21})& (a_{i2}, b_{22})&\cdots & (a_{ik}, b_{2k})\cr
\vdots & \vdots & \vdots & \vdots\cr
(a_{i1}, b_{n1})& (a_{i2}, b_{n2})&\cdots & (a_{ik}, b_{nk})}.
\end{eqnarray*}

\begin{theorem} \label{Comp-Thm}
Suppose that $A$
is a super-simple
$\mathrm{MCA}_{d_1+1}(N_1; t, k, (v_1,\break  v_2, \ldots, v_k))$ of optimum
size $N_1=(d_1+1) \prod_{i=k-t+1}^k v_i$
and $B$ is a super-simple
$\mathrm{MCA}_{d_2+1}(N_2; t, k, (u_1, u_2, \ldots, u_k))$ of optimum
size $N_2=\break (d_2+1) \prod_{i=k-t+1}^k u_i$.
Then $A\otimes B$ is a $(d,t)$-$\mathrm{DTA}(N_1N_2; k, (v_1u_1, v_2u_2, \ldots,\break  v_ku_k))$ of optimum
size $N_1N_2$ attaining the lower bound in (\ref{Eqn-LowB-M}), where
$d=(d_1+1)(d_2+1)-1$.
\end{theorem}

\begin{corollary} \label{Comp-Cor3}
Suppose that there exists a super-simple
$\mathrm{MCA}_{d+1}(N; t,\break  k, (v_1, v_2, \ldots, v_k))$ of optimum
size $N=(d+1) \prod_{i=k-t+1}^k v_i$. Then:
\begin{longlist}[(2)]
\item[(1)] there exists a
$(\lambda(d+1)-1,t)$-$\mathrm{DTA}(\lambda m^t N; k, (mv_1, mv_2, \ldots,
mv_k))$ of optimum
size $\lambda m^t N$ meeting the lower bound in (\ref{Eqn-LowB-M}) for
any integer $\lambda$ with
$2\le\lambda\le m$, provided that an $\mathrm{OA}(t+1, k+1, m)$ exists;

\item[(2)] there exists a $(d,t)$-$\mathrm{DTA}(m^t N; k, (mv_1, mv_2, \ldots,mv_k))$ of optimum
size $m^t N$ meeting the lower bound in (\ref{Eqn-LowB-M}),
provided that an $\mathrm{OA}(t, k, m)$ exists.
\end{longlist}
\end{corollary}

Theorem \ref{Comp-Thm} and Corollary \ref{Comp-Cor3}
provide the basic idea for constructing
new optimum detecting arrays meeting the lower bound
in (\ref{Eqn-LowB-M}) from the old ones.
Combining with the existence results for
the ingredient arrays, we obtain the following theorem.

\begin{theorem} \label{Comp-Ex}
There exists:
\begin{longlist}[(3)]
\item[(1)] a $(2^m-1,2)$-$\mathrm{DTA}(18^m;4,(2^m)^1(3^m)^3)$ achieving the
lower bound in (\ref{Eqn-LowB-M}) for any integer $m \geq2$;\vadjust{\goodbreak}
\item[(2)] a $(2\lambda-1,2)$-$\mathrm{DTA}(18\lambda v^2;4,(2v)^1(3v)^3)$
achieving the lower bound\break  in~(\ref{Eqn-LowB-M}) for any integer
$\lambda\leq v$, where $v\ge4$ and $v\not\equiv2$ ($\mathrm{mod}$ 4); and
\item[(3)] a $(1,2)$-$\mathrm{DTA}(18v^2;4,(2v)^1(3v)^3)$ achieving the lower
bound in (\ref{Eqn-LowB-M}) for any integer $v\neq2, 4, 6$.
\end{longlist}
\end{theorem}

\subsection{Existence of optimum DTAs with $k=t+1$}

We begin with the following known results.

\begin{lemma}\label{known-1}
 Suppose that $2\le v_1\le v_2\le v_3
\leq v_4$. Then there
exists:
\begin{longlist}[(2)]
\item[(1)] an $\mathrm{MCA}(N; 2,3,(v_1, v_2,v_3))$ of optimum size $N=v_2v_3$
[\citet{MSSW2003}];
\item[(2)] an $\mathrm{MCA}(N; 3,4,(v_1, v_2,v_3, v_4))$ of optimum size
$N=v_2v_3 v_4$ [\citet{CSWY2011}].
\end{longlist}
\end{lemma}

Taking advantage of some composition constructions, we
can completely determine the existence spectrum for
a $(d, t)$-$\mathrm{DTA}(N; t+1, (v_1, v_2, \ldots, v_{t+1}))$ of
optimum size $N$ meeting the
lower bound in (\ref{Eqn-LowB-M}).

\begin{lemma}\label{main-con1}
Suppose that an $\mathrm{MCA}(N;t,k, (v_1, v_2, \ldots, v_i, \ldots, v_k))$ and
an $\mathrm{MCA}(N';t-1,k-1, (v_1, v_2, \ldots, v_{i-1}, v_{i+1}, \ldots, v_k))$
all exist, where \mbox{$1 \le i\le k$}.
Then an $\mathrm{MCA}(N'',t,k, (v_1, v_2, \ldots, v_{i-1}, (v_i+e), v_{i+1},
\ldots, v_k))$
exists for any positive integer $e$, where $N''=N+eN'$.
\end{lemma}

\begin{lemma}\label{main-con2}
Suppose that $2\leq v_1\leq v_2\leq\cdots\leq v_{t+1}$.
Then an optimum $\mathrm{MCA}(N; t,t+1, (v_1, v_2, \ldots, v_{t+1}))$ with
$N=\prod_{i=2}^{t+1}v_i$ exists
for any integer $t\ge2$.
\end{lemma}

\begin{theorem}\label{main}
Suppose\vspace*{2pt} that $2\leq v_1\leq v_2\leq\cdots\leq v_{t+1}$. Then an
optimum $(d-1,t)$-$\mathrm{DTA}(N; t+1,(v_1,v_2,\ldots,v_{t+1}))$ with $N=\prod
_{i=2}^{t+1}v_i$ exists if and only if $d\leq
v_1$.
\end{theorem}

\section{Conclusion}\label{s6}

Detecting arrays are a special class of covering arrays.
However, detecting arrays, especially mixed level detecting arrays,
have much more complicated structure than that of usual covering arrays.
They are introduced
to generate test suites that are capable of
identifying and determining the faulty interactions between factors
in a component-based system. Compared to using covering arrays,
this is a more useful approach
in an application context.
In this paper, the notion of ``$d$-extendible''
is proposed and used to
explore the combinatorial features of mixed level detecting arrays.
By way of equivalent combinatorial configurations, an optimality criterion
is found.
As a consequence, combinatorial approaches and heuristic algorithms
are employed to
produce a number of examples and infinite classes of mixed level
detecting arrays, which
are optimum in size meeting the lower bound in (\ref{Eqn-LowB-M}).

\begin{appendix}\label{app}
\section*{Appendix: Proofs}

\begin{pf*}{Proof of Theorem \ref{Thm-LowB-M}}
Let $A$ be a $(d,t)$-$\mathrm{DTA}(N; k, (v_1, v_2, \ldots, v_k))$. From
Lemma \ref{lbound-d>1},
we know that $|\rho(A,T)|\geq d+1 $
for any $t$-way interaction $T$.
Since there are $v_{_{k-t+1}}\times v_{_{k-t+2}}\times\cdots
\times v_k$ different $t$-way interactions,
$\{(r, x_r)\dvtx  k-t+1\leq r\leq k \}$,
$A$ must contain $(d + 1) \prod_{i=k-t+1}^k v_i$ rows, which implies
$(d,t)$-$\mathrm{DTAN}(k, (v_1, v_2, \ldots, v_k)) \ge(d+1)\prod_{i=k-t+1}^k
v_i$.
\end{pf*}

\begin{pf*}{Proof of Lemma \ref{TheOAEqu}}
``$\Rightarrow$'' If $A$ is not super-simple, then
there exists a $(t+1)$-way interaction $T_{t+1,1}$ such that
$|\rho(A,T_{t+1,1})| \ge2$. Select the first $t$ elements of
$T_{t+1,1}$ to form a $t$-way interaction $T_t$, then $T_{t+1,1}$ is
an extension of $T_t$. Let { ${\mathcal T} = \{T| T \mbox{ is an
extension of } T_t \mbox{ with } |\rho(A,T)|\ge1\}$}.
Obviously, $T_{t+1,1} \in{\mathcal T}$. Notice that $|\rho
(A,T_{t+1,1})| \ge2$
and $|\rho(A,T_{t})| = d+1$ as $A$ is an $\mathrm{OA}_{d+1}(t,k,v)$,
we have $|{\mathcal T}| \le d$. However,
$\rho(A,T_{t}) \setminus\rho(A,{\mathcal T})
=\varnothing$. This is a contradiction to the $d$-{extendible} property.

``$\Leftarrow$'' Suppose there exist
a $t$-way interaction $T_t$ and its $d$
extensions, $A_{t+1,i}$, $i=1,\ldots, d$,
such that $\rho(A,T_{t}) \setminus\bigcup_{i=1} ^d \rho(A,T_{t+1,i})
=\varnothing$. Since $A$ is an $\mathrm{OA}_{d+1}(t,k,v)$, $|\rho(A,T_{t})| =
d+1$, there must exist $i \in\{1,2,\ldots,d\}$, such that $|\rho
(A,T_{t+1,i})| \ge2$.
This is a contradiction to the {super-simple property.}
\end{pf*}

\begin{pf*}{Proof of Theorem \ref{TheEquivd}}
``$\Leftarrow$'' If $A$ is not a $(d,t)$-$\mathrm{DTA}(N; k, (v_1, v_2,
\ldots,\break  v_k))$, then there exist a $t$-way interaction $T$ and
a set of $d$ $t$-way interactions $\mathcal{T}=\{T_1,T_2,\ldots, T_d\}$
such that $T\notin
\mathcal{T}$ and $\rho(A,T) \subseteq
\rho(A,\mathcal{T})$.
If there exists $T_i \in\mathcal{T}$,
such that all first coordinates of elements of $T_i$
and $T$ coincide, then define any extension of $T$
as $T_{t+1,i}$.
On the other hand, for each $T_i \in\mathcal{T}$,
if there exists an element of $T_i$, say $(j_i,\sigma_{j_i})$, such
that the first coordinate of any element of $T$ is not~$j_i$,
then add the element $(j_i,\sigma_{j_i})$ to $T$ to form an extension
of $T$,
denoted by $T_{t+1,i}$.
Now for any row $R \in\rho(A,T) \subseteq
\rho(A,\mathcal{T})$, if $R \in\rho(A,T_i)$,
then $R \in\rho(A,T_{t+1,i})$, which tells
$\rho(A,T_{t}) \setminus\bigcup_{i=1} ^d \rho(A,T_{t+1,i})
=\varnothing$. Notice that for $i\neq j$,
$T_{{t+1,i}}$ and $T_{{t+1,j}}$ can be the same $(t+1)$-way interaction.

``$\Rightarrow$'' Otherwise, there exist a $t$-way
interaction $T$ and $d$ extensions of $T$, denoted by $T_{t+1,i}$,
$1\leq i\leq d$,
such that $\rho(A,T_{t}) \setminus\bigcup_{i=1} ^d
\rho(A,T_{t+1,i})=\varnothing$. Then $\rho(A,T_{t}) = \bigcup_{i=1}
^d \rho(A,T_{t+1,i})$,
as $\rho(A,T_{t+1,i}) \subseteq\rho(A,T)$ for $1\leq i \leq d$ by
Lemma~\ref{T,T+1}. For each
extension $T_{t+1,i}$, we can select a $t$-way interaction
$T_{t,i}\neq T$ such that $T_{t+1,i}$ is their common extension.
Notice that for $i\neq j$, $T_{t,i}$ and $T_{t,j}$ cannot be the same
$t$-way interaction,
otherwise, $T_{t+1,i}$ and $T_{t+1,j}$ must be the same $(t+1)$-way interaction.
Now we have
$\rho(A,T) = \bigcup_{i=1} ^d \rho(A,T_{t+1,i}) \subseteq
\bigcup_{i=1} ^d \rho(A,T_{t,i})$. This is { a contradiction}
to the definition of a $(d,t)$-detecting array.
\end{pf*}

\begin{pf*}{Proof of Theorem \ref{Efeature1}}
``$\Rightarrow$'' If $A$ is an {optimum}
$(d,t)$-$\mathrm{DTA}(N; t+1,(v_1, v_2, \ldots, v_{t+1}))$, then $A$
is ``$d$-extendible.'' By Lemma \ref{lbound-d>1}, $|\rho(A,T)|\geq
d+1$ for
any $t$-way interaction $T$. Thus, $A$ is an {optimum}
$\mathrm{MCA}_{d+1}(N;t, t+1, (v_1, v_2, \ldots, v_{t+1}))$. In the
following, we prove that $A$ is super-simple. Otherwise, suppose
that there are {two identical row vectors}, denoted by
$R_1=(x_1,x_2,\ldots, x_{t+1})$. Let $T=(x_2,x_3,\ldots,x_{t+1})$.
$T$ occurs at least $d+1$ times as the subarray indexed by the columns
$\{2,3,\ldots, t+1\}$, and hence exactly
$d+1$ times [since $A$ contains precisely $(d+1)\prod_{i=2}^{t+1} v_i$
rows]. Write $T_{t+1,i}$ $(i=1,2,3,\ldots, d)$ for
the $d$ extensions of $T$, {where $T_{t+1,1}=R_1$.} Then $\rho(A,T)
\setminus\bigcup_{i=1}^d \rho(A,T_{t+1,i}) =\varnothing$,
a contradiction to $A$ being ``$d$-extendible.''

``$\Leftarrow$'' Let $B$ be a super-simple $\mathrm{MCA}_{d+1}(N;t, t+1, (v_1,
v_2, \ldots, v_{t+1}))$ with
$N=(d+1)\prod_{i=2}^{t+1} v_i$. We prove that $B$ is a $d$-{extendible}
$\mathrm{MCA}_{d+1}(N;t, k, (v_1,\break  v_2, \ldots,v_k))$. If not, there exists a
$t$-way interaction $T$ and $d$ extensions $T_{t+1,i}$ $(i=1,2,3,\ldots, d)$
for $T$ such that $\rho(A,T) \setminus\bigcup_{i=1}^d \rho
(A,T_{t+1,i}) =\varnothing$.
Since $T$ occurs at least $d+1$ times,
there is one $T_{t+1,i} $ such that
$|\rho(A,T_{t+1,i})|\geq2$. This implies that there is at least
one $(t + 1)$-tuple of symbols occurring in some $t + 1$ columns as
rows more than once. This is a contradiction to the
super-simple property of $B$. The conclusion then follows from Theorem
\ref{TheEquivd}.
\end{pf*}

\begin{pf*}{Proof of Theorem \ref{Efeature2}}
The proof can be divided into
two cases, one is similar to that of Theorem \ref{Efeature1},
and the other is actually a copy of the proof of Theorem 2.4 in \citet{STY2011}.
\end{pf*}

\begin{pf*}{Proof of Theorem \ref{Feature}}
The proof is similar to the counterpart in Theorem~\ref{Efeature1}.
\end{pf*}

\begin{pf*}{Proof of Lemma \ref{Lem-Fixed}}
\citet{TY2011} { proved} that
an optimum $(1,2)$-$\mathrm{DTA}(N; k, q)$ with $N = 2q^2$
is equivalent to a super-simple orthogonal array
$\mathrm{OA}_{2}(2,k,q)$. So for any 2-way interaction
$T$, $|\rho(A,T)| = 2$. Moreover, the pair of rows
covering $T$ must be different for {each distinct 2-way}
interaction $T$. The number of 2-way interactions
for a $k$-factor array is ${k \choose2} \times q^2 = k(k-1)q^2/2$.
{However,} there are totally
${N \choose2} = q^2(2q^2-1)$ pairs of rows. Thus, we
have $k (k-1) q^2/2 \le q^2 (2q^2-1)$, which implies
$k \le2q$ as $k$ and $q$ are both positive integers.
\end{pf*}

\begin{pf*}{Proof of Lemma \ref{Lem-Existence}}
(i) We first prove that either $u$ or $k$ is less than $2$. Otherwise,
if both $u$ and $k$ are greater than or equal to $2$, then we have
$N=6w$. Let $A_0$ be the subarray with the last column taking the value
$0$. Then the number of rows of $A_0$ is 6. For
three-level factors of $A_0$, there are
two rows taking level zero, one and two; for
two-level factors of $A_0$, there are
three rows taking level zero and one.
Moreover, each column of the factors with the same level of $A_0$ cannot
be a level permutation of another one. Thus,
two two-level factors occupy $ 2\times2\times{{3}\choose{2}} -2 =
10$ pairs
of rows, on which there is at least one two-level factor taking
the same level. Two three-level factors occupy $2 \times3 = 6$ pairs
of rows of~$A_0$, on which there is at least one three-level factor taking
the same level. These $10+6=16$ pairs of
rows are mutually
different, which is a contradiction to
the ${{6}\choose{2}} = 15$ possible pairs of six rows of $A_0$.

(ii) If $k=0$, then $N=4w$. Let $A_0$ be the subarray with the last
column taking the value $0$.
Then the number of rows of $A_0$ is 4. For columns of $A_0$ other than
the last one, there are
two zeros and two ones. Moreover, each of these zeros or ones in the
same column
occupies a
different pair of rows of $A_0$,
otherwise we can always find a two-way interaction $T$ (which contains
the last column as a component) and one of its
extension $T_a$, so that the rows covering $T$ and $T_a$
are identical.
Thus, we must have ${{4}\choose{2}} \ge2 u$, that is, $u \le3$.

(iii) If $k=1$, then $N=6w$. Let $A_0$ be the subarray with the last
column taking the value $0$.
Then the number of rows of $A_0$ is 6. Without loss of generality, for
the three-level factor of $A_0$, assume the first two rows
take level zero, the middle two take level one, and
the last two take level two.
For columns of $A_0$ other than the last two,
the first two, the middle and the last two
rows must take distinct levels, that is, one takes zero and the other
takes one.
The number of such possible columns is $2 \times2 \times2 = 8$.
Finally, among all these $u$ columns, one cannot
be a level permutation of another. Thus,
$u \le4$.

(iv)
If $u=0$ and $k \ge1$, then $N=6w$.
Let $A_0$ be the subarray with the last column taking the value $0$.
Then the number of rows of $A_0$ is 6. For columns of $A_0$ other than
the last one, there are
two zeros, two ones and two twos. Moreover,
each distinct element in the same column of $A_0$ {occupies a}
different pair of rows of $A_0$. Thus, we must have ${{6}\choose{2}}
\ge3 k$, that is, $k \le5$.

(v)
If $u=1$ and $k \ge1$, then $N=6w$.
Let $A_0$ be the subarray with the last column taking the value $0$.
Then the number of rows of $A_0$ is 6.
The two-level factor occupies $ 2\times{{3}\choose{2}} = 6$ pairs
of rows, on which the two-level factor takes
the same level. Thus, $k\times3 + 6 \le{{6}\choose{2}}$,
which gives $k \le3$.
\end{pf*}

\begin{pf*}{Proof of Theorem \ref{Comp-Thm}}
Write $A=(a_{ij})$ $(i \in I_{N_1}, j\in I_k)$
whose entries on the $j$th column are taken from $\Z_{v_j}$ for $1\le
j\le k$.
Write $B=(b_{ij})$ $(i \in I_{N_2}, j\in I_k)$
whose entries on the $j$th column are taken from $\Z_{u_j}$ for $1\le
j\le k$.
Then, by definition, the Kronecker product $C$ is {an} $(N_1N_2)\times
k$ matrix
whose entries on the $j$th column are taken from the Cartesian product
$\Z_{v_j}\times\Z_{u_j}$ for $1\le j\le k$. Furthermore, each row
$(a_{i1}, a_{i2}, \ldots, a_{ik})$ of
$A$ generates $N_2$ rows in the Kronecker product $C$ whose projections
on the second
component are precisely the $N_2$ rows of $B$.
So the rows of each $(N_1N_2)\times t$ subarray of
$C$ cover all $t$-tuples
of values from the $t$ columns at least $(d_1+1)(d_2+1)$ times, as $A$
and $B$ are both MCAs of strength $t$
with indices $d_1+1$ and $d_2+1$,
respectively.
Similarly, the super-simple property of $A$ and $B$
guarantees that $C$ is super-simple. This is because the projection on
the first
component of each row in any $t+1$ columns
of $C$ is a row of the
corresponding columns of $A$, while the projection on the second
component is a row of the corresponding
columns of $B$. It follows that $C$ is a super-simple
$\mathrm{MCA}_{d+1}(N_1N_2; k, (v_1u_1, v_2u_2, \ldots, v_ku_k))$, where the size
$N_1N_2=(d_1+1)(d_2+1) \prod_{i=k-t+1}^k (v_iu_i)$. The conclusion
then follows from
Theorem 3.12.
\end{pf*}

\begin{pf*}{Proof of Corollary \ref{Comp-Cor3}}
From \citet{STY2011}, we know that an $\mathrm{OA}(t + 1, k + 1, m)$ implies
the existence of a super-simple $\mathrm{OA}_{\lambda}(t, k, m)$ with
$2\le\lambda\le m$.
{So the} conclusion (1) holds
by taking $u_1=u_2=\cdots=u_k=m$ in Theorem \ref{Comp-Thm}.
The conclusion (2) follows directly
from Theorem \ref{Comp-Thm} by taking $u_1=u_2=\cdots=u_k=m$ and $d_2=0$.
\end{pf*}

\begin{pf*}{Proof of Corollary \ref{Comp-Ex}}
(1) Apply Theorem \ref{Comp-Thm} with a super-simple $\mathrm{MCA}_2(18; 2,4,2^13^3)$
given in Section~\ref{sec4};
(2) Apply the conclusion $(1)$ in Corollary \ref{Comp-Cor3}
with a super-simple $\mathrm{MCA}_2(18;2,4,2^13^3)$ and an $\mathrm{OA}(3,5,v)$ given in
\citet{JY2010};
(3) Apply the conclusion $(2)$ in Corollary \ref{Comp-Cor3} with a
super-simple $\mathrm{MCA}_2(18;2,4,2^13^3)$ and
an $\mathrm{OA}(2,4,v)$ given in \citet{CD2007} and \citet{HSS1999}.
\end{pf*}

\begin{pf*}{Proof of Lemma \ref{main-con1}}
Let $A$ be an $\mathrm{MCA}(N;t,k, (v_1, v_2, \ldots, v_i, \ldots, v_k))$
whose entries on the $j$th column are taken from $\Z_{v_j}$ for $1\le
j\le k$.
For any particular value $i$ with $1\le i\le k$, we form an
$\mathrm{MCA}(N';t-1,k, (v_1, v_2, \ldots,\break v_{i-1}, 1, v_{i+1}, \ldots,
v_k))$ by
simply inserting one constant column
\[
(v_i-1+r, v_i-1+r, \ldots, v_i-1+r)^T
\]
in the $i$th position of the given $\mathrm{MCA}(N';t-1,k-1, (v_1, v_2, \ldots,
v_{i-1}, v_{i+1}, \ldots,\break  v_k))$.
For given positive integer $e$, we do this for $r=1, 2, \ldots, e$.
This produces $e$ MCAs of index unity whose values on the $i$th column
are taken from $\{v_i, v_i+1, \ldots, v_i+e-1\}$.
We then concatenate the obtained $e$ MCAs together with the MCA $A$ to form
an $\mathrm{MCA}(N'',t,k, (v_1, v_2, \ldots, v_{i-1}, (v_i+e), v_{i+1}, \ldots, v_k))$
with $N''=N+eN'$, as desired.
\end{pf*}

\begin{pf*}{Proof of Lemma \ref{main-con2}}
The proof is done by mathematical induction on $t$. For $t=2,3$, the
conclusion follows from Lemma \ref{known-1}.
{Assume that} the conclusion holds
when $t =n\ge3$ and consider the case $t=n+1$.
It is well known that {both an $\mathrm{OA}(n+1,n+2,v_1)$ and an $\mathrm{OA}(n,n+1,v_1)$}
exist for any given integer $v_1$ [\citet{HSS1999}].
In essence, they are an optimum CA$(v_1^{n+1}; n+1,n+2, v_1)$ and
an optimum CA$(v_1^n; n,n+1, v_1)$, respectively. When $v_2>v_1$, we
apply Lemma \ref{main-con1}
with $i=2,e=v_2-v_1,(t,k)=(n+1,n+2)$ and these two MCAs. This creates
an $\mathrm{MCA}(N_2;n+1,n+2,(v_1,v_2,v_1,\ldots,v_1))$
of optimum size,
denoted by $A_2$, where $N_2=v_1^{n+1}+(v_2-v_1)v_1^{n}=v_1^{n}v_2$. In
the case $v_1=v_2$,
$A_2$ can be simply taken to be an $\mathrm{OA}(n+1,n+2,v_1)$. By the induction
hypothesis, we have also
an $\mathrm{MCA}(v_1^{n}v_2; n,n+1, (v_1, v_2, \ldots, v_{n+1}))$ denoted by $B$.
{So we} can again apply Lemma \ref{main-con1}, as above,
with $(t,k)=(n+1,n+2)$, $i=3, e=v_3-v_1$ and the MCAs $A_2$, $B$ to
form an $\mathrm{MCA}(N_3;n+1,n+2,(v_1,v_2,v_3,v_1,\ldots,v_1))$
of optimum size $N_3=v_1^{n}v_2+(v_3-v_1)v_1^{n-1}v_2=v_1^{n-1}v_2v_3$,
denoted by $A_3$.
This process can be continued recursively until an optimum
$\mathrm{MCA}(N_{n+2}; n+1,n+2, (v_1, v_2, \ldots, v_{n+2}))$ is obtained,
where $N_{n+2}=\prod_{i=2}^{n+2}v_i$. Consequently, the assertion also
holds when $t=n+1$. This completes the proof.
\end{pf*}

\begin{pf*}{Proof of Theorem \ref{main}}
By taking advantage of the equivalence described in Theorem \ref
{Efeature1}, we need\vspace*{1pt}
only to show that there exists a super-simple
$\mathrm{MCA}_{d}(d\cdot N; t, t+1, (v_1, v_2, \ldots, v_{t+1}))$ of optimum
size $d\cdot N=d\cdot{\prod_{i=2}^{t+1} v_i}$.
From Lemma~\ref{main-con2}, we know that an
$\mathrm{MCA}(N;t,t+1,(v_1,v_2,\ldots,v_{t+1}))$
of optimum size $N={\prod_{i=2}^{t+1} v_i}$ exists. Write $A$ for such
an MCA.
For each $i$ with $0\le i\le d-1$, we form a new optimum
$\mathrm{MCA}(N;t,t+1,(v_1,v_2,\ldots,v_{t+1}))$, denoted by $A_{\sigma^i}$, by
permuting the entries in
the first column of $A$ with
the permutation $\sigma^i$. Write $D=(A_{\sigma^{0}}^T\cdots
A_{\sigma^{d-1}}^T)^T$
for the\vspace*{2pt} concatenation of the obtained $d$ MCAs. It is obvious that $D$
is an $\mathrm{MCA}_d(dN;t,t+1,(v_1,v_2,\ldots,v_{t+1}))$.
It is {optimum}, since its size $d\cdot N=d\cdot{\prod_{i=2}^{t+1} v_i}$.
It is now left to show that $D$ is
super-simple. In fact, otherwise, there would be a $(t+1)$-tuple
of values from certain $t+1$ columns of $D$
occurring as rows of these columns at least twice. Without loss of
generality, we assume
that the $(t+1)$-tuple $\underline{a}=(a_1,a_2,\ldots, a_{t+1})$ of
values from the first $t+1$ columns of $D$
occurs twice, as rows of these columns.
This implies that there must be $i$ and $j$ with $0\le i\ne j \le d-1$
so that the
$(t+1)$-tuple $\underline{a}$ occurs in the corresponding columns of
both $A_{\sigma^i}$ and $A_{\sigma^j}$,
since each MCA $A_{\sigma^r}$ $(0\le r\le d-1)$ has index $\lambda=1$.
Let $b_1$ be an arbitrary value from the first
column of $D$. Notice that the MCA $A$ is of strength $t$ and index
$\lambda=1$.
By considering the $t+1$-tuple $\underline{b}=(b_1,a_2,\ldots, a_{t+1})$,
one can see that $\sigma^i(b_1)=a_1=\sigma^j(b_1)$ from the
construction of $A_{\sigma^i}$ and $A_{\sigma^j}$.
This leads $\sigma^i=\sigma^j$. It is a contradiction to the
definition of $\sigma$ (as the order of $\sigma$ is $v_1$).
\end{pf*}
\end{appendix}

\section*{Acknowledgments}
The authors would like to thank the Associate Editor and two referees
for their comments and suggestions, which were greatly beneficial for
this paper.


%

\printaddresses

\begin{thebibliography}{26}

\bibitem[\protect\citeauthoryear{Brownlie, Plowse and Phadke}{1992}]{BPP1992}
%
\begin{barticle}[auto:STB|2014/05/26|13:19:10]
\bauthor{\bsnm{Brownlie},~\bfnm{J.}\binits{J.}},
\bauthor{\bsnm{Plowse},~\bfnm{J.}\binits{J.}} \AND
\bauthor{\bsnm{Phadke},~\bfnm{M.}\binits{M.}}
(\byear{1992}).
\btitle{Robust Testing of AT\&T PMX/StarMail using OATS}.
\bjournal{ATT Tech. J.}
\bvolume{3}
\bpages{41--47}.
\end{barticle}
%
\bptok{imsref}%
\endbibitem

%
%

\bibitem[\protect\citeauthoryear{Cohen et~al.}{1996}]{CDPP1996}
%
\begin{barticle}[auto:STB|2014/05/26|13:19:10]
\bauthor{\bsnm{Cohen},~\bfnm{D.~M.}\binits{D.~M.}},
\bauthor{\bsnm{Dalal},~\bfnm{S.~R.}\binits{S.~R.}},
\bauthor{\bsnm{Parelius},~\bfnm{J.}\binits{J.}} \AND
\bauthor{\bsnm{Patton},~\bfnm{G.~C.}\binits{G.~C.}}
(\byear{1996}).
\btitle{The combinatorial design approach to automatic test generation}.
\bjournal{IEEE Softw.}
\bvolume{13}
\bpages{83--88}.
\end{barticle}
%
\bptok{imsref}\vadjust{\goodbreak}%
\endbibitem

\bibitem[\protect\citeauthoryear{Colbourn}{2004}]{Col2004}
%
\begin{barticle}[mr]
\bauthor{\bsnm{Colbourn},~\bfnm{Charles~J.}\binits{C.~J.}}
(\byear{2004}).
\btitle{Combinatorial aspects of covering arrays}.
\bjournal{Matematiche (Catania)}
\bvolume{59}
\bpages{125--172}.
\bid{mr={2243029}}
\end{barticle}
%
\bptok{imsref}%
\endbibitem

\bibitem[\protect\citeauthoryear{Colbourn and Dinitz}{2007}]{CD2007}
%
\begin{bbook}[auto]
\bauthor{\bsnm{Colbourn},~\bfnm{C. J.}\binits{C. J.}} \AND
\bauthor{\bsnm{Dinitz},~\bfnm{J. H.}\binits{J. H.}}
(\byear{2007}).
\btitle{The CRC Handbook of Combinatorial Designs}.
\bpublisher{Chapman \& Hall/CRC},
\blocation{Boca Raton, FL}.
\end{bbook}
%
\bptok{imsref}%
\endbibitem

\bibitem[\protect\citeauthoryear{Colbourn and McClary}{2008}]{CM2008}
%
\begin{barticle}[mr]
\bauthor{\bsnm{Colbourn},~\bfnm{Charles~J.}\binits{C.~J.}} \AND
\bauthor{\bsnm{McClary},~\bfnm{Daniel~W.}\binits{D.~W.}}
(\byear{2008}).
\btitle{Locating and detecting arrays for interaction faults}.
\bjournal{J. Comb. Optim.}
\bvolume{15}
\bpages{17--48}.
\bid{doi={10.1007/s10878-007-9082-4}, issn={1382-6905}, mr={2375213}}
\end{barticle}
%
\bptok{imsref}%
\endbibitem

\bibitem[\protect\citeauthoryear{Colbourn et~al.}{2006}]{CMMSSY2006}
%
\begin{barticle}[mr]
\bauthor{\bsnm{Colbourn},~\bfnm{Charles~J.}\binits{C.~J.}},
\bauthor{\bsnm{Martirosyan},~\bfnm{Sosina~S.}\binits{S.~S.}},
\bauthor{\bsnm{Mullen},~\bfnm{Gary~L.}\binits{G.~L.}},
\bauthor{\bsnm{Shasha},~\bfnm{Dennis}\binits{D.}},
\bauthor{\bsnm{Sherwood},~\bfnm{George~B.}\binits{G.~B.}} \AND
\bauthor{\bsnm{Yucas},~\bfnm{Joseph~L.}\binits{J.~L.}}
(\byear{2006}).
\btitle{Products of mixed covering arrays of strength two}.
\bjournal{J. Combin. Des.}
\bvolume{14}
\bpages{124--138}.
\bid{doi={10.1002/jcd.20065}, issn={1063-8539}, mr={2202133}}
\end{barticle}
%
\bptok{imsref}%
\endbibitem

\bibitem[\protect\citeauthoryear{Colbourn et~al.}{2011}]{CSWY2011}
%
\begin{barticle}[mr]
\bauthor{\bsnm{Colbourn},~\bfnm{Charles~J.}\binits{C.~J.}},
\bauthor{\bsnm{Shi},~\bfnm{Ce}\binits{C.}},
\bauthor{\bsnm{Wang},~\bfnm{Chengmin}\binits{C.}} \AND
\bauthor{\bsnm{Yan},~\bfnm{Jie}\binits{J.}}
(\byear{2011}).
\btitle{Mixed covering arrays of strength three with few factors}.
\bjournal{J. Statist. Plann. Inference}
\bvolume{141}
\bpages{3640--3647}.
\bid{doi={10.1016/j.jspi.2011.05.018}, issn={0378-3758}, mr={2817370}}
\end{barticle}
%
\bptok{imsref}%
\endbibitem

\bibitem[\protect\citeauthoryear{Dalal et~al.}{1999}]{DKLPH}
%
\begin{bincollection}[auto:STB|2014/05/26|13:19:10]
\bauthor{\bsnm{Dalal},~\bfnm{S.~R.}\binits{S.~R.}},
\bauthor{\bsnm{Karunanithi},~\bfnm{A.~J.~N.}\binits{A.~J.~N.}},
\bauthor{\bsnm{Leaton},~\bfnm{J.~M.~L.}\binits{J.~M.~L.}},
\bauthor{\bsnm{Patton},~\bfnm{G.~C.~P.}\binits{G.~C.~P.}} \AND
\bauthor{\bsnm{Horowitz},~\bfnm{B.~M.}\binits{B.~M.}}
(\byear{1999}).
\btitle{Model-based testing in practice}.
In \bbooktitle{Proceedings of the International Conference on Software
Engineering (ICSE'99)}
\bpages{285--294}.
\bpublisher{ACM}, \blocation{New York}.
\end{bincollection}
%
\bptok{imsref}%
\endbibitem

\bibitem[\protect\citeauthoryear{Esmailzadeh et~al.}{2011}]{ETMJ2011}
%
\begin{barticle}[mr]
\bauthor{\bsnm{Esmailzadeh},~\bfnm{Nabaz}\binits{N.}},
\bauthor{\bsnm{Talebi},~\bfnm{Hooshang}\binits{H.}},
\bauthor{\bsnm{Momihara},~\bfnm{Koji}\binits{K.}} \AND
\bauthor{\bsnm{Jimbo},~\bfnm{Masakazu}\binits{M.}}
(\byear{2011}).
\btitle{A new series of main effects plus one plan for {$2^m$}
factorial experiments with {$m=4\lambda\pm1$} and {$2m$} runs}.
\bjournal{J. Statist. Plann. Inference}
\bvolume{141}
\bpages{1567--1574}.
\bid{doi={10.1016/j.jspi.2010.11.008}, issn={0378-3758}, mr={2747925}}
\end{barticle}
%
\bptok{imsref}%
\endbibitem

\bibitem[\protect\citeauthoryear{Ghosh and Burns}{2001}]{GB2001}
%
\begin{barticle}[mr]
\bauthor{\bsnm{Ghosh},~\bfnm{Subir}\binits{S.}} \AND
\bauthor{\bsnm{Burns},~\bfnm{Colleen}\binits{C.}}
(\byear{2001}).
\btitle{Two general classes of search designs for factor screening
experiments with factors at three levels}.
\bjournal{Metrika}
\bvolume{54}
\bpages{1--17}.
\bid{doi={10.1007/PL00003989}, issn={0026-1335}, mr={1863415}}
\end{barticle}
%
\bptok{imsref}%
\endbibitem

\bibitem[\protect\citeauthoryear{Hedayat, Sloane and Stufken}{1999}]{HSS1999}
%
\begin{bbook}[mr]
\bauthor{\bsnm{Hedayat},~\bfnm{A.~S.}\binits{A.~S.}},
\bauthor{\bsnm{Sloane},~\bfnm{N.~J.~A.}\binits{N.~J.~A.}} \AND
\bauthor{\bsnm{Stufken},~\bfnm{John}\binits{J.}}
(\byear{1999}).
\btitle{Orthogonal Arrays}.
\bpublisher{Springer},
\blocation{New York}.
\bid{doi={10.1007/978-1-4612-1478-6}, mr={1693498}}
\end{bbook}
%
\bptok{imsref}%
\endbibitem

\bibitem[\protect\citeauthoryear{Ji and Yin}{2010}]{JY2010}
%
\begin{barticle}[mr]
\bauthor{\bsnm{Ji},~\bfnm{Lijun}\binits{L.}} \AND
\bauthor{\bsnm{Yin},~\bfnm{Jianxing}\binits{J.}}
(\byear{2010}).
\btitle{Constructions of new orthogonal arrays and covering arrays of
strength three}.
\bjournal{J. Combin. Theory Ser. A}
\bvolume{117}
\bpages{236--247}.
\bid{doi={10.1016/j.jcta.2009.06.002}, issn={0097-3165}, mr={2592899}}
\end{barticle}
%
\bptok{imsref}%
\endbibitem

\bibitem[\protect\citeauthoryear{Kuhn, Wallace and Gallo}{2004}]{KWG}
%
\begin{barticle}[auto:STB|2014/05/26|13:19:10]
\bauthor{\bsnm{Kuhn},~\bfnm{D.~R.}\binits{D.~R.}},
\bauthor{\bsnm{Wallace},~\bfnm{D.~R.}\binits{D.~R.}} \AND
\bauthor{\bsnm{Gallo},~\bfnm{A.~M.}\binits{A.~M.}}
(\byear{2004}).
\btitle{Software fault interactions and implications for software testing}.
\bjournal{IEEE Trans. Softw. Eng.}
\bvolume{30}
\bpages{418--421}.
\end{barticle}
%
\bptok{imsref}%
\endbibitem

%
%

\bibitem[\protect\citeauthoryear{Moura et~al.}{2003}]{MSSW2003}
%
\begin{barticle}[mr]
\bauthor{\bsnm{Moura},~\bfnm{Lucia}\binits{L.}},
\bauthor{\bsnm{Stardom},~\bfnm{John}\binits{J.}},
\bauthor{\bsnm{Stevens},~\bfnm{Brett}\binits{B.}} \AND
\bauthor{\bsnm{Williams},~\bfnm{Alan}\binits{A.}}
(\byear{2003}).
\btitle{Covering arrays with mixed alphabet sizes}.
\bjournal{J. Combin. Des.}
\bvolume{11}
\bpages{413--432}.
\bid{doi={10.1002/jcd.10059}, issn={1063-8539}, mr={2012427}}
\end{barticle}
%
\bptok{imsref}%
\endbibitem

\bibitem[\protect\citeauthoryear{Shi, Tang and Yin}{2012}]{STY2011}
%
\begin{barticle}[mr]
\bauthor{\bsnm{Shi},~\bfnm{Ce}\binits{C.}},
\bauthor{\bsnm{Tang},~\bfnm{Yu}\binits{Y.}} \AND
\bauthor{\bsnm{Yin},~\bfnm{Jianxing}\binits{J.}}
(\byear{2012}).
\btitle{The equivalence between optimal detecting arrays and
super-simple {OA}s}.
\bjournal{Des. Codes Cryptogr.}
\bvolume{62}
\bpages{131--142}.
\bid{doi={10.1007/s10623-011-9498-9}, issn={0925-1022}, mr={2886266}}
\end{barticle}
%
\bptok{imsref}%
\endbibitem

\bibitem[\protect\citeauthoryear{Shirakura, Takahashi and
Srivastava}{1996}]{STS1996}
%
\begin{barticle}[mr]
\bauthor{\bsnm{Shirakura},~\bfnm{Teruhiro}\binits{T.}},
\bauthor{\bsnm{Takahashi},~\bfnm{Tadashi}\binits{T.}} \AND
\bauthor{\bsnm{Srivastava},~\bfnm{Jagdish~N.}\binits{J.~N.}}
(\byear{1996}).
\btitle{Searching probabilities for nonzero effects in search designs
for the noisy case}.
\bjournal{Ann. Statist.}
\bvolume{24}
\bpages{2560--2568}.
\bid{doi={10.1214/aos/1032181169}, issn={0090-5364}, mr={1425968}}
\end{barticle}
%
\bptok{imsref}%
\endbibitem


\bibitem[\protect\citeauthoryear{Srivastava}{1975}]{Srivastava1975}
%
\begin{bincollection}[mr]
\bauthor{\bsnm{Srivastava},~\bfnm{J.~N.}\binits{J.~N.}}
(\byear{1975}).
\btitle{Designs for searching non-negligible effects}.
In \bbooktitle{A Survey of Statistical Design and Linear Models
({P}roc. {I}nternat. {S}ympos., {C}olorado {S}tate {U}niv., {F}t.
{C}ollins, {C}olo., 1973)}
\bpages{507--519}.
\bpublisher{North-Holland},
\blocation{Amsterdam}.
\bid{mr={0398023}}
\end{bincollection}
%
\bptok{imsref}%
\endbibitem

\bibitem[\protect\citeauthoryear{Tang and Yin}{2011}]{TY2011}
%
\begin{barticle}[mr]
\bauthor{\bsnm{Tang},~\bfnm{Yu}\binits{Y.}} \AND
\bauthor{\bsnm{Yin},~\bfnm{Jian~Xing}\binits{J.~X.}}
(\byear{2011}).
\btitle{Detecting arrays and their optimality}.
\bjournal{Acta Math. Sin. (Engl. Ser.)}
\bvolume{27}
\bpages{2309--2318}.
\bid{doi={10.1007/s10114-011-0184-7}, issn={1439-8516}, mr={2853789}}
\end{barticle}
%
\bptok{imsref}%
\endbibitem

\end{thebibliography}
\end{document}